\documentclass[10pt]{amsart}

\usepackage{tikz}
\usepackage{amssymb,amsmath,latexsym,graphicx}
\usepackage{charter,eucal}
\usepackage{amssymb}
\usepackage[all,cmtip]{xy}
\usepackage{rotating,multicol}
\usetikzlibrary{decorations.markings}

\newtheorem{theorem}{Theorem }[section]

\newtheorem{lemma}[theorem]{Lemma}
\newtheorem{observation}[theorem]{Observation}

\newtheorem{remark}[theorem]{Remark}
\newtheorem{corollary}[theorem]{Corollary}
\newtheorem{proposition}[theorem]{Proposition}
\newtheorem{principle}[theorem]{\textsc{Principle}}

\newcommand{\bt}{\begin{theorem}}
\newcommand{\et}{\end{theorem}}
\newcommand{\bmt}{\begin{maintheorem}}
\newcommand{\emt}{\end{maintheorem}}
\newcommand{\bc}{\begin{corollary}}
\newcommand{\bl}{\begin{lemma}}
\newcommand{\ec}{\end{corollary}}
\newcommand{\el}{\end{lemma}}
\newcommand{\bo}{\begin{observation}}
\newcommand{\eo}{\end{observation}}
\newcommand{\bp}{\begin{proposition}}
\newcommand{\ep}{\end{proposition}}
\newcommand{\br}{\begin{remark}}
\newcommand{\er}{\end{remark}}
\newcommand{\bpr}{\begin{principle}}
\newcommand{\epr}{\end{principle}}

\def\PG{\mathbf{PG}}

\def\eop{\hspace*{\fill}$\blacksquare$}

\usepackage{amssymb} 
\usepackage{amsmath} 
\usepackage[utf8]{inputenc} 
\usepackage{latexsym}
\title {\bf Regular pseudo-hyperovals and regular pseudo-ovals in even characteristic}
\author {J. A. THAS}
\address{Ghent University, Department of Mathematics, Krijgslaan 281, S22, B-9000 Ghent, Belgium}
\email{jat@cage.ugent.be}
\begin {document}
\maketitle
\begin {abstract}
S. Rottey and G. Van de Voorde characterized regular pseudo-ovals of $\PG(3n - 1, q)$, $q = 2^h$, $h > 1$ and $n$ prime. Here an alternative proof is given and slightly stronger results are obtained.
\end{abstract}

\section{Introduction}
Pseudo-ovals and pseudo-hyperovals were introduced in \cite{JAT: 71}; see also \cite{TTVM: 06}. These objects play a key role in the theory of translation generalized quadrangles \cite{PT: 09, TTVM: 06}. Pseudo-hyperovals only exist in even characteristic. A characterization of regular pseudo-ovals in odd characteristic was given in \cite{CTW: 85}; see also \cite{TTVM: 06}. In \cite{RVDV: 15} a characterization of regular pseudo-ovals and regular pseudo-hyperovals in $\PG(3n - 1, q)$, $q$ even, $q \neq 2$ and $n$ prime, is obtained. Here a shorter proof is given and slightly stronger results are obtained.

\section{Ovals and hyperovals}
A {\em $k$-arc} in $\PG(2, q)$ is a set of $k$ points, $k \ge 3$, no three of which are collinear. Any non-singular conic of $\PG(2, q)$ is a $(q + 1)$-arc. If $\mathcal {K}$ is any $k$-arc of $\PG(2, q)$, then $k \leq {q + 2}$. For $q$ odd $k \leq q+1$ and for $q$ even a $(q + 1)$-arc extends to a $(q + 2)$-arc; see \cite{JWPH: 98}. A $(q + 1)$-arc is an {\em oval}; a $(q + 2)$-arc, $q$ even, is a {\em complete oval} or {\em hyperoval}.

A famous theorem of B. Segre \cite{S: 54} tells us that for $q$ odd every oval of $\PG(2, q)$ is a non-singular conic. For $q$ even, there are many ovals that are not conics \cite{JWPH: 98}; also, there are many hyperovals that do not contain a conic \cite{JWPH: 98}.

\section{Generalized ovals and hyperovals}
Arcs, ovals and hyperovals can be generalized by replacing their points with $m$-dimensional subspaces to obtain generalized $k$-arcs, generalized ovals and generalized hyperovals. These have strong connections to generalized quadrangles, projective planes, circle geometries, flocks and other structures. See \cite{PT: 09, TTVM: 06, JAT: 71, JAT: 11, CTW: 85, PVDV: 13}. Below, some basic definitions and results are formulated; for an extensive study, many applications and open problems, see \cite{TTVM: 06}.

A {\em generalized $k$-arc} of $\PG(3n - 1, q)$, $n \geq 1$, is a set of $k$ $(n - 1)$-dimensional subspaces of $\PG(3n - 1, q)$ every three of which generate $\PG(3n - 1, q)$. If $q$ is odd then $k \leq q^n + 1$, if $q$ is even then $k \leq {q^n + 2}$. Every generalized $(q^n + 1)$-arc of $\PG(3n - 1, q)$, $q$ even, can be extended to a generalized $(q^n + 2)$-arc.

If $\mathcal {O}$ is a generalized $(q^n + 1)$-arc in $\PG(3n - 1, q)$, then it is a  {\em pseudo-oval} or {\em generalized oval} or \lbrack $n - 1$\rbrack-{\em oval} of $\PG(3n - 1, q)$. For $n = 1$, a [0]-oval is just an oval of $\PG(2, q)$. If $\mathcal{O}$ is a generalized $(q^n + 2)$-arc in $\PG(3n - 1, q)$, $q$ even, then it is a {\em pseudo-hyperoval} or {\em generalized hyperoval} or \lbrack $n - 1$\rbrack-{\em hyperoval} of $\PG(3n - 1, q)$. For $n = 1$, a [0]-hyperoval is just a hyperoval of $\PG(2, q)$.

If $\mathcal{O} = \lbrace \pi_0, \pi_1, \cdots, \pi_{q^n} \rbrace$ is a pseudo-oval of $\PG(3n - 1, q)$, then $\pi_i$ is contained in exactly one $(2n - 1)$-dimensional subspace $\tau_i$ of $\PG(3n - 1, q)$ which has no point in common with $(\pi_0\cup\pi_1\cup \cdots \cup\pi_{q^n})\backslash\pi_i$, with $i = 0, 1, \cdots, q^n$; the space $\tau_i$ is the {\em tangent space} of $\mathcal{O}$ at $\pi_i$. For $q$ even the $q^n + 1$ tangent spaces of $\mathcal{O}$ contain a common $(n - 1)$-dimensional space $\pi_{q^n + 1}$, the {\em nucleus} of $\mathcal{O}$; also, $\mathcal{O} \cup \lbrace\pi_{q^n + 1}\rbrace$  is a pseudo-hyperoval of $\PG(3n - 1, q)$. For $q$ odd, the tangent spaces of a pseudo-oval $\mathcal{O}$ are the elements of a pseudo-oval $\mathcal{O^\ast}$ in the dual space of $\PG(3n - 1, q)$.

\section{Regular pseudo-ovals and pseudo-hyperovals}
In the extension $\PG(3n - 1, q^n)$ of $\PG(3n - 1, q)$, consider $n$ subplanes $\xi_i$, $i = 1, 2, \cdots, n$, that are conjugate in the extension $\mathbb{F}_{q^n}$ of $\mathbb{F}_q$ and which span $\PG(3n - 1, q^n)$. This means that they form an orbit of the Galois group corresponding to this extension and span $\PG(3n - 1, q^n)$.

In $\xi_1$ consider an oval $\mathcal{O}_1 = \lbrace x_0^{(1)}, x_1^{(1)}, \cdots, x_{q^n}^{(1)}\rbrace$. Further, let $x_i^{(1)}, x_i^{(2)}, \cdots, x_i^{(n)}$, with $i = 0, 1, \cdots, q^n$, be conjugate in $\mathbb{F}_{q^n}$ over $\mathbb{F}_q$. The points $x_i^{(1)}, x_i^{(2)}, \cdots, x_i^{(n)}$ define an $(n - 1)$-dimensional subspace $\pi_i$ over $\mathbb{F}_q$ for $i = 0, 1, \cdots, q^n$. Then, $\mathcal{O} = \lbrace \pi_0, \pi_1, \cdots, \pi_{q^n}\rbrace$ is a generalized oval of $\PG(3n - 1, q)$. These objects are the {\em regular} or {\em elementary pseudo-ovals}. If $\mathcal{O}_1$ is replaced by a hyperoval, and so $q$ is even, then the corresponding $\mathcal{O}$ is a {\em regular} or {\em elementary pseudo-hyperoval}.

All known pseudo-ovals and pseudo-hyperovals are regular.

\section{Characterizations}
Let $\mathcal{O} = \lbrace \pi_0, \pi_1, \cdots, \pi_{q^n} \rbrace$ be a pseudo-oval in $\PG(3n - 1, q)$. The tangent space of $\mathcal{O}$ at $\pi_i$ will be denoted by $\tau_i$, with $i = 0, 1, \cdots, q^n$. Choose $\pi_i$, $i \in \lbrace 0, 1, \cdots, q^n \rbrace$, and let $\PG(2n - 1, q) \subseteq \PG(3n - 1, q)$ be skew to $\pi_i$. Further, let $\tau_i \cap \PG(2n - 1, q) = \eta_i$ and $\langle  \pi_i, \pi_j \rangle \cap \PG(2n - 1, q) = \eta_j$, with $j \neq i$. Then $\lbrace \eta_0, \eta_1, \cdots, \eta_{q^n} \rbrace = \Delta_i$ is an $(n - 1)$-spread of $\PG(2n - 1, q)$.

Now, let $q$ be even and let $\pi$ be the nucleus of $\mathcal{O}$. Let $\PG(2n - 1, q) \subseteq \PG(3n - 1, q)$ be skew to $\pi$. If $\zeta_j = \PG(2n - 1, q) \cap \langle \pi, \pi_j \rangle$, then $\lbrace \zeta_0, \zeta_1, \cdots, \zeta_{q^n} \rbrace = \Delta$ is an $(n - 1)$-spread of $\PG(2n - 1, q)$.

Next, let $q$ be odd. Choose $\tau_i$, with $i \in \lbrace 0, 1, \cdots, q^n \rbrace$. If $\tau_i \cap \tau_j = \delta_j$, with $j \neq i$, then $\lbrace \delta_0, \delta_1, \cdots, \delta_{i - 1}, \pi_i, \delta_{i + 1}, \cdots , \delta _{q^n} \rbrace = \Delta^\star_i$ is an $(n - 1)$-spread of $\tau_i$.

Finally, let $q$ be even and let $\mathcal{O} = \lbrace \pi_0, \pi_1, \cdots, \pi_{q^n + 1} \rbrace$ be a pseudo-hyperoval in $\PG(3n - 1, q)$. Choose $\pi_i$, with $i \in \lbrace 0, 1, \cdots, q^n + 1 \rbrace$, and let $\PG(2n - 1, q) \subseteq \PG(3n - 1, q)$ be skew to $\pi_i$. Let $\langle \pi_i, \pi_j \rangle \cap \PG(2n - 1, q) = \eta_j$, with $j \neq i$. Then $\lbrace \eta_0, \eta_1, \cdots, \eta_{i - 1}, \eta_{i + 1}, \cdots,  \eta_{q^n + 1} \rbrace = \Delta_i$ is an $(n - 1)$- spread of $\PG(2n - 1, q)$.

\begin{theorem}[Casse, Thas and Wild \cite{CTW: 85}]
Consider a pseudo-oval $\mathcal{O}$ with $q$ odd. Then at least one of the $(n - 1)$-spreads $\Delta_0, \Delta_1, \cdots, \Delta_{q^n}, \Delta^\star_0, \Delta^\star_1, \cdots, \Delta^\star_{q^n}$ is regular if and only if they all are regular if and only if the pseudo-oval $\mathcal{O}$ is regular. In such a case $\mathcal{O}$ is essentially a conic over $\mathbb{F}_{q^n}$.
\end{theorem}

\begin{theorem}[Rottey and Van de Voorde \cite{RVDV: 15}]
Consider a pseudo-oval $\mathcal{O}$ in $\PG(3n - 1, q)$ with $q = 2^h$, $h > 1$, $n$ prime. Then $\mathcal{O}$ is regular if and only if all $(n - 1)$-spreads $\Delta_0, \Delta_1, \cdots, \Delta_{q^n}$ are regular.
\end{theorem}

\section{Alternative proof and improvements}
\begin{theorem}
Consider a pseudo-hyperoval $\mathcal{O}$ in $\PG(3n - 1, q)$, $q = 2^h, h > 1$ and $n$ prime. Then $\mathcal{O}$ is regular if and only if all $(n - 1)$-spreads $\Delta_i$, with $i = 0, 1, \cdots, q^n + 1$, are regular.
\end{theorem}

{\em Proof}. \quad
If $\mathcal{O}$ is regular, then clearly all $(n - 1)$-spreads $\Delta_i$, with $i = 0, 1, \cdots, q^n + 1$, are regular.

Conversely, assume that the $(n - 1)$-spreads $\Delta_0, \Delta_1, \cdots, \Delta_{q^n + 1}$ are regular. Let $\mathcal{O} = \lbrace \pi_0, \pi_1, \cdots, \pi_{q^n + 1} \rbrace$ and let $\mathcal{\hat{O}} = \lbrace \beta_0, \beta_1, \cdots, \beta_{q^n + 1} \rbrace$ be the dual of $\mathcal{O}$, with $\beta_i$  being the dual of $\pi_i$.

Choose $\beta_i, i \in \lbrace 0, 1, \cdots, q^n + 1 \rbrace$, and let $\beta_i \cap \beta_j = \alpha_{ij}, j \neq i$. Then 
\begin{equation}
 \lbrace \alpha_{i0}, \alpha_{i1}, \cdots, \alpha_{i, i - 1}, \alpha_{i, i + 1}, \cdots, \alpha_{i, q^n + 1} \rbrace = \Gamma_i 
 \end{equation}
 is an $(n - 1)$-spread of $\beta_i$.

Now consider $\beta_i, \beta_j, \Gamma_i, \Gamma_j, \alpha_{ij}, j \neq i $. In $\Gamma_j$ we next consider a $(n - 1)$-regulus $\gamma_j$ containing $\alpha_{ij}$. The $(n - 1)$-regulus $\gamma_j$ is a set of maximal spaces of a Segre variety $\mathcal{S}_{1; n - 1}$; see Section 4.5 in \cite{HT: 16}. The $(n - 1)$-regulus $\gamma_j$ and the $(n - 1)$-spread $\Gamma_i$ of $\beta_i$ generate a regular $(n - 1)$-spread $\Sigma (\gamma_j, \Gamma_i)$ of $\PG(3n - 1, q)$. This can be seen as follows. The elements of $\Gamma_i$ intersect $n$ lines $U_1, U_2, \cdots, U_n$ which are conjugate in $\mathbb F_{q^n}$ over $\mathbb F_q$, that is, they form an orbit of the Galois group corresponding to this extension. Let $\alpha_{ij} \cap U_l = \lbrace u_l \rbrace$, with $l = 1, 2, \cdots, n$. Now consider the transversals $T_1, T_2, \cdots, T_n$ of the elements of $\gamma_j$, with $T_l$ containing $u_l$. The $n$ planes $T_lU_l = \theta_l$ intersect all elements of $\gamma_j$ and $\Gamma_i$. The $(n - 1)$-dimensional subspaces of $\PG(3n - 1, q)$ intersecting $\theta_1, \theta_2, \cdots, \theta_n$ are the elements of the regular $(n - 1)$-spread $\Sigma (\gamma_j, \Gamma_i)$. The elements of this spread are the points of a plane $\PG(2, q^n)$, with as lines the $(2n - 1)$-dimensional spaces containing at least two (and then $q^n + 1$) elements of the spread. Hence the $q + 2$ elements of $\mathcal{\hat{O}}$ containing an element of $\gamma_j$, say $\beta_i = \beta_{i_1}, \beta_{i_2}, \cdots, \beta_{i_{q + 1}}, \beta_{i_{q + 2}} = \beta_j$, are lines of $\PG(2, q^n)$. Dualizing, the elements $\pi_{i_1}, \pi_{i_2}, \cdots, \pi_{i_{q + 2}}$ are points of $\PG(2, q^n)$.

Now consider $\beta_{i_2}$ and $\gamma_j$, and repeat the argument. Then there arise $n$ planes $\theta'_l$ intersecting all elements of $\gamma_j$ and $\Gamma_{i_2}$. The $(n - 1)$-dimensional subspaces of $\PG(3n - 1, q)$ intersecting $\theta'_1, \theta'_2, \cdots, \theta'_n$ are the elements of the regular $(n - 1)$-spread $\Sigma(\gamma_j, \Gamma_{i_2})$. The elements of this spread are the points of a plane $\PG'(2, q^n)$, with as lines the $(2n - 1)$-dimensional spaces containing $q^n + 1$ elements of the spread. Hence $\beta_{i_1}, \beta_{i_2}, \cdots, \beta_{i_{q + 2}}$ are lines of $\PG'(2, q^n)$. Dualizing, the elements $\pi_{i_1}. \pi_{i_2}, \cdots, \pi_{i_{q + 2}}$ are points of $\PG'(2, q^n)$.

First, assume that $\lbrace \theta_1, \theta_2, \cdots, \theta_n \rbrace \cap \lbrace \theta'_1, \theta'_2, \cdots, \theta'_n \rbrace = \emptyset$. Consider $\pi_{i_1}, \pi_{i_2}, \pi_{i_3}, \pi_{i_4}$. The planes of $\PG(3n - 1, q^n)$ intersecting these four spaces constitute a set $\mathcal{M}$ of maximal spaces of a Segre variety $\mathcal{S}_{2; n - 1}$ \cite{B: 61}. The planes $\theta_1, \theta_2, \cdots, \theta_n, \theta'_1, \theta'_2, \cdots, \theta'_n$ are elements of $\mathcal{M}$. It follows that $(\theta_1 \cup \theta_2 \cup \cdots \cup \theta_n) \cap (\theta'_1 \cup \theta'_2 \cup \cdots \cup \theta'_n) = \emptyset$.


Consider any $(n - 1)$-dimensional subspace $\pi \in \lbrace \pi_{i_5}, \pi_{i_6}, \cdots, \pi_{i_{q + 2}} \rbrace$ of $\PG(3n - 1, q)$. We will show that $\pi$ is a maximal subspace of $\mathcal{S}_{2;n - 1}$. Let $\theta_i \cap \pi_j = \lbrace t_{ij} \rbrace, \theta_i^\prime  \cap \pi_j = \lbrace t_{ij}^\prime \rbrace, i = 1, 2, \cdots, n, j= i_1, i_2, \cdots, i_{q+2}$. If $t_{ij_1}t_{ij_2} \cap t_{ij_3}t_{ij_4} = \lbrace v_i \rbrace, t_{ij_1}^\prime t_{ij_2}^\prime \cap t_{ij_3}^\prime t_{ij_4}^\prime = \lbrace v_i^\prime \rbrace$, with $j_1, j_2, j_3, j_4$ distinct, then $v_1, v_2, \cdots, v_n$ are conjugate and similarly $v_1^\prime, v_2^\prime, \cdots, v_n^\prime$ are conjugate. Hence $\langle v_1, v_2, \cdots, v_n \rangle = \langle v_1^\prime, v_2^\prime, \cdots, v_n^\prime \rangle$ defines a $(n - 1)$-dimensional space over $\mathbb{F}_q$ which intersects $\theta_1, \theta_2, \cdots, \theta_n^\prime$ (over $\mathbb{F}_{q^n}$).The points $t_{ij}$, with $j = i_1, i_2, \cdots, i_{q + 2}$, generate a subplane of $\theta_i$, and the points $t_{ij}^\prime$, with $j = i_1, i_2, \cdots, i_{q + 2}$, generate a subplane of $\theta_i^\prime$, with $i = 1, 2, \cdots, n$. Let $q = 2^h$ and let $\mathbb{F}_{2^v}$ be the subfield of $\mathbb{F}_{q^n} = \mathbb{F}_{2^{hn}} $ over which these subplanes are defined; so $v \vert hn$. Then $v < hn$ as otherwise the spreads of $\PG(3n - 1, q)$ defined by $\theta_1, \theta_2, \cdots, \theta_n$ and $\theta_1^\prime, \theta_2^\prime, \cdots, \theta_n^\prime$ coincide, clearly not possible. The $(n - 1)$-regulus $\gamma_j$ implies that the subplanes contain a line over $\mathbb{F}_q$, so $h \vert v$. As $n$ is prime we have $v =h$, so $2^v = q$. Hence the $2n$ subplanes are defined over $\mathbb{F}_q$. It follows that the $q + 2$ elements $\pi_{i_1}, \pi_{i_2}, \cdots, \pi_{i_{q + 2}}$ are maximal subspaces of the Segre variety $\mathcal{S}_{2; n - 1}$. Hence $\pi$ is a maximal subspace of $\mathcal{S}_{2; n -1}$. It follows that $\pi_1, \pi_2, \cdots, \pi_{q + 2}$ are maximal subspaces of $\mathcal{S}_{2; n - 1}$.

Now consider a $\PG(2, q )$ which intersects $\pi_{i_1}, \pi_{i_2}, \pi_{i_3}, \pi_{i_4}$. The $(n - 1)$-dimensional spaces $\pi_{i_1}, \pi_{i_2}, \cdots, \pi_{i_{q + 2}}$ are maximal spaces of $\mathcal{S}_{2; n - 1}$ which intersect $\PG(2, q)$; they are maximal spaces of the Segre variety $\mathcal{S}_{2; n - 1} \cap \PG(3n - 1,q)$ of $\PG(3n - 1, q)$.

Consider $\pi_{i_1}$ and also a $\PG(2n - 1, q)$ skew to $\pi_{i_1}$. If we project $\pi_{i_2}, \pi_{i_3}, \cdots, \pi_{i_{q + 2}}$ from $\pi_{i_1}$ onto $\PG(2n - 1, q)$, then by the foregoing paragraph the $q + 1$ projections constitute a $(n - 1)$-regulus of $\PG(2n - 1, q)$. Similarly, if we project from $\pi_{i_s}$, $s$ any element of $\lbrace 1, 2, \cdots, q + 2 \rbrace$. Equivalently, if $s \in \lbrace 1, 2, \cdots, q + 2 \rbrace$ then the spaces $\beta_{i_s} \cap \beta_{i_t}$, with $t = 1, 2, \cdots, s - 1, s + 1, \cdots, q + 2$, form a $(n - 1)$-regulus of $\beta_{i_s}$.

Now assume that the condition $\lbrace \theta_1, \theta_2, \cdots, \theta_n \rbrace \cap \lbrace \theta'_1, \theta'_2, \cdots, \theta'_n \rbrace = \emptyset$ is satisfied for any choice of $\beta_i, \beta_j, \gamma_j, \beta_{i_2}$. In such a case every $(n - 1)$-regulus contained in a spread $\Gamma_s$ defines a Segre variety $\mathcal{S}_{2; n - 1}$ over $\mathbb{F}_q$. Let us define the following design $\mathcal{D}$. Points of $\mathcal{D}$ are the elements of $\mathcal{\hat{O}}$, a block of $\mathcal{D}$ is a set of $q + 2$ elements of $\mathcal{\hat{O}}$, containing at least one space of a $(n - 1)$-regulus contained in some regular spread $\Gamma_s$, and incidence is containment. Then $\mathcal{D}$ is a $4 - (q^n + 2, q + 2, 1)$ design. By Kantor \cite{K: 74} this implies that $q = 2$, a contradiction.

Consequently, we may assume that for at least one quadruple $\beta_i, \beta_j, \gamma_j, \beta_{i_2}$ we have
\begin{equation}
 \lbrace \theta_1, \theta_2, \cdots, \theta_n \rbrace = \lbrace \theta'_1, \theta'_2, \cdots, \theta'_n \rbrace. 
 \end{equation}
 In such a case the $q ^n + 2$ elements of $\mathcal{\hat{O}}$ are lines of the plane $\PG(2, q^n)$. It follows that $\mathcal{O}$ is regular. \eop
 
 \begin{theorem}
 Consider a pseudo-oval $\mathcal{O}$ in $\PG(3n - 1, q)$, with $q = 2^h, h > 1$ and $n$ prime. Then $\mathcal{O}$ is regular if and only if all $(n - 1)$-spreads $\Delta_0, \Delta_1, \cdots, \Delta_{q^n}$ are regular.
 \end{theorem}
 
 {\em Proof}. \quad
 If $\mathcal{O}$ is regular, then clearly all $(n - 1)$-spreads $\Delta_0, \Delta_1, \cdots, \Delta_{q^n}$ are regular. 
 
 Conversely, assume that the $(n - 1)$-spreads $\Delta_0, \Delta_1, \cdots, \Delta_{q^n}$ are regular. Let $\mathcal{O} = \lbrace \pi_0, \pi_1, \cdots, \pi_{q^n} \rbrace$, let $\pi_{q^n + 1}$ be the nucleus of $\mathcal{O}$, let $\mathcal{\bar{O}} = \mathcal{O} \cup \lbrace \pi_{q^n + 1} \rbrace$, let $\mathcal{\hat{O}}$ be the dual of $\mathcal{O}$, let $\mathcal{\hat{\bar{O}}}$ be the dual of $\mathcal{\bar{O}}$, and let $\beta_i$ be the dual of $\pi_i$.
 
 Choose $\beta_i, i \in \lbrace 0, 1, \cdots, q^n + 1 \rbrace$, and let $\beta_i \cap \beta_j = \alpha_{ij}, j \neq i$. Then 
 \begin{equation}
 \lbrace \alpha_{i0}, \alpha_{i1}, \cdots, \alpha_{i, i - 1}, \alpha_{i, i + 1}, \cdots, \alpha_{i, q^n + 1} \rbrace = \Gamma_i
 \end{equation}
 is an $(n - 1)$-spread of $\beta_i$.
 
 Now consider $\beta_i, \beta_j, \Gamma_i, \Gamma_j, \alpha_{ij}$, with $j \neq i$ and $i, j \in \lbrace 0, 1, \cdots, q^n \rbrace$. In $\Gamma_j$  we next consider a $(n - 1)$-regulus $\gamma_j$ containing $\alpha_{ij}$ and $\alpha_{j, q^n + 1}$. The $(n - 1)$-regulus $\gamma_j$ is a set of maximal spaces of a Segre variety $\mathcal{S}_{1; n - 1}$. The $(n - 1)$-regulus $\gamma_j$ and the $(n - 1)$-spread $\Gamma_i$ of $\beta_i$ generate a regular $(n - 1)$-spread $\Sigma (\gamma_j, \Gamma_i)$ of $\PG(3n - 1, q)$. Such as in the proof of Theorem 6.1 we introduce the elements $U_l, u_l, T_l, \theta_l, l = 1, 2, \cdots, n$, and the plane $\PG(2, q^n)$. The $q + 2$ elements of $\mathcal{\hat{\bar{O}}}$ containing an element of $\gamma_j$, say $\beta_i = \beta_{i_1}, \beta_{i_2}, \cdots, \beta_{i_q}, \beta_j = \beta_{i_{q + 1}}, \beta_{q^n + 1}$, are lines of $\PG(2, q^n)$. Dualizing, the elements $\pi_{i_1}, \pi_{i_2}, \cdots, \pi_{i_{q + 1}}, \pi_{q^n + 1}$ are points of $\PG(2, q^n)$.
 
 Now consider $\beta_{i_2}$ and $\gamma_j$, and repeat the argument. Then there arise $n$ planes $\theta'_l$ of $\PG(3n - 1, q^n)$ intersecting all elements of $\gamma_j$ and $\Gamma_{i_2}$, and a $(n - 1)$-spread $\Sigma(\gamma_j, \Gamma_{i_2})$ of $\PG(3n - 1, q)$. The elements of this spread are the points of a plane $\PG'(2, q^n)$. The spaces $\beta_{i_1}, \beta_{i_2}, \cdots, \beta_{i_{q + 1}}, \beta_{q^n + 1}$ are lines of $\PG'(2, q^n)$. Dualizing, the elements $\pi_{i_1}, \pi_{i_2}, \cdots, \pi_{i_{q + 1}}, \pi_{q^n + 1}$ are points of $\PG'(2, q^n)$.
 
 First, assume that $\lbrace \theta_1, \theta_2, \cdots, \theta_n \rbrace \cap \lbrace \theta'_1, \theta'_2, \cdots, \theta'_n \rbrace = \emptyset$. Consider $\pi_{i_1}, \pi_{i_2}, \pi_{i_3}, \pi_{i_4}$. The planes of $\PG(3n - 1, q^n)$ intersecting these four spaces constitute a set $\mathcal{M}$ of maximal spaces of a Segre variety $\mathcal{S}_{2; n - 1}$. The planes $\theta_1, \theta_2, \cdots, \theta_n, \theta'_1, \theta'_2, \cdots, \theta'_n$ are elements of $\mathcal{M}$. It follows that $(\theta_1 \cup \theta_2 \cup \cdots \cup \theta_n) \cap (\theta'_1 \cup \theta'_2 \cup \cdots \cup \theta'_n) = \emptyset$. Let $\pi \in \lbrace \pi_{i_5}, \pi_{i_6}, \cdots, \pi_{i_{q + 1}}, \pi_{q^n + 1} \rbrace$. As in the proof of Theorem 6.1 one shows that $\pi$ is a maximal subspace of $\mathcal{S}_{2; n - 1}$. It follows that $\pi_{i_1}, \pi_{i_2}, \cdots, \pi_{i_{q + 1}}, \pi_{q^n + 1}$ are maximal subspaces of $\mathcal{S}_{2; n - 1}$.
 
Next consider a $\PG(2, q)$ which intersects $\pi_{i_1}, \pi_{i_2}, \pi_{i_3}, \pi_{i_4}$. The $(n - 1)$-dimensional spaces $\pi_{i_1}, \pi_{i_2}, \cdots, \pi_{i_{q + 1}}, \pi_{q^n + 1}$ are maximal spaces of $\mathcal{S}_{2; n - 1}$ which intersect the plane $\PG(2, q)$; they are maximal spaces of the Segre variety $\mathcal{S}_{2; n - 1} \cap \PG(3n - 1, q)$ of $\PG(3n - 1, q)$. Such as in the proof of Theorem 6.1 it follows that the spaces $\beta_{q^n + 1} \cap \beta_{i_t}$, with $t = 1, 2, \cdots, q + 1$, form a $(n - 1)$-regulus of $\beta_{q^n + 1}$. 
 
 Now assume that the condition $\lbrace \theta_1, \theta_2, \cdots, \theta_n \rbrace \cap \lbrace \theta'_1, \theta'_2, \cdots, \theta'_n \rbrace = \emptyset$ is satisfied for any choice of $\beta_i, \beta_j, \gamma_j, \beta_{i_2}$, $j \neq i$ and $i, j \in \lbrace 0, 1, \cdots, q^n \rbrace$. Let $\alpha_1, \alpha_2, \alpha_3$ be distinct elements of $\Gamma_{q^n + 1}$. Then $\beta_i, \beta_j, \gamma_j, \beta_{i_2}$ can be chosen in such a way that $\alpha_1 \in \beta_i, \alpha_2 \in \beta_j, \alpha_2 \in \gamma_j, \beta_{i_2} \cap \beta_j \in \gamma_j$ with $\alpha_3 \in \beta_{i_2}$. Hence the $(n - 1)$-regulus in $\beta_{q^n + 1}$ defined by $\alpha_1, \alpha_2, \alpha_3$ is subset of $\Gamma_{q^n + 1}$. From \cite{HT: 16} now follows that the $(n - 1)$-spread $\Gamma_{q^n + 1}$ of $\beta_{q^n + 1}$ is regular. By Theorem 6.1 the pseudo-hyperoval $\mathcal{\bar{O}}$ is regular, and so $\mathcal{O}$ is regular. But in such a case the condition $\lbrace \theta_1, \theta_2, \cdots, \theta_n \rbrace \cap \lbrace \theta'_1, \theta'_2, \cdots, \theta'_n \rbrace = \emptyset$ is never satisfied, a contradiction.
 
Consequently, we may assume that for at least one quadruple $\beta_i, \beta_j, \gamma_j, \beta_{i_2}$ we have $\lbrace \theta_1, \theta_2, \cdots, \theta_n \rbrace = \lbrace \theta'_1, \theta'_2, \cdots, \theta'_n \rbrace$. In such a case the $q^n + 2$ elements of $\mathcal{\hat{\bar{O}}}$ are lines of the plane $\PG(2, q^n)$. It follows that $\mathcal{\bar{O}}$, and hence also $\mathcal{O}$, is regular. \eop
 
 \begin{theorem}
 Consider a pseudo-hyperoval $\mathcal{O}$ in $\PG(3n - 1, q)$, $q = 2^h, h > 1$ and $n$ prime. Then $\mathcal{O}$ is regular if and only if at least $q^n - 1$ elements of $\lbrace \Delta_0, \Delta_1, \cdots, \Delta_{q^n + 1} \rbrace$ are regular.
 \end{theorem}
 
 {\em Proof}. \quad
 If $\mathcal{O}$ is regular, then clearly all $(n - 1)$-spreads $\Delta_i$, with $i = 0, 1, \cdots, q^n + 1$, are regular.
 
 Conversely, assume that $\rho$, with $\rho \geq q^n - 1$, elements of $\lbrace \Delta_0, \Delta_1, \cdots, \Delta_{q^n + 1} \rbrace$ are regular.
 
 If $\rho = q^n + 2$, then $\mathcal{O}$ is regular by Theorem 6.1; if $\rho = q^n + 1$, then $\mathcal{O}$ is regular by Theorem 6.2.
 
 Now assume that $\rho = q^n$ and that $\Delta_2, \Delta_3, \cdots, \Delta_{q^n + 1}$ are regular. We have to prove that $\Delta_0$ is regular. We use the arguments in the proof of Theorem 6.2. If one of the elements $\alpha_1, \alpha_2, \alpha_3$, say $\alpha_1$, in the proof of Theorem 6.2 is $\beta_0 \cap \beta_1$, then let $\gamma_j$ contain $\beta_j \cap \beta_i, \beta_j \cap \beta_0, \beta_j \cap \beta_1$ and let $\beta_{i_2} \neq \beta_1$, with $i, j \in \lbrace 2, 3, \cdots, q^n + 1 \rbrace$. Now see the proof of the preceding theorem.
 
 Finally, assume that $\rho = q^n - 1$ and that $\Delta_3, \Delta_4, \cdots, \Delta_{q^n + 1}$ are regular. We have to prove that $\Delta_0$ is regular. We use the arguments in the proof of Theorem 6.2. If exactly one of the elements $\alpha_1, \alpha_2, \alpha_3$, say $\alpha_1$, in the proof of Theorem 6.2 is $\beta_0 \cap \beta_1$ or $\beta_0 \cap \beta_2$, then proceed as in the preceding paragraph with $\beta_{i_2} \neq \beta_1, \beta_2$. Now assume that two of the elements $\alpha_1, \alpha_2, \alpha_3$, say $\alpha_1$ and $\alpha_2$, are $\beta_0 \cap \beta_1$ and $\beta_0 \cap \beta_2$. Now consider all $(n - 1)$-reguli in $\Delta_0$ containing $\alpha_1$ and $\alpha_3$, and assume, by way of contradiction, that no one of these $(n - 1)$-reguli contains $\alpha_2$. The number of these $(n - 1)$-reguli is $\frac{q^n - 2} {q - 1}$, and so $q = 2$, a contradiction. It follows that the $(n - 1)$-regulus in $\beta_0$ defined by $\alpha_1, \alpha_2, \alpha_3$ is contained in $\Delta_0$. Now we proceed as in the proof of Theorem 6.2. \eop
 
 \section{Final remarks}
 7.1. {\bf {The cases $q = 2$ and $n$ not prime}} \\
 For $q = 2$ or $n$ not prime other arguments have to be developed.
 
 7.2. {\bf {Improvement of Theorem 6.3}} \\
 Let $\mathcal{D} = (P, B, \in)$ be an incidence structure satisfying the following conditions.
 \begin{itemize}
 \item[{\rm (i)}]
 $\vert P \vert = q^n + 1$, $q$ even, $q \neq 2$;
 \item[{\rm (ii)}]
 the elements of $B$ are subsets of size $q + 1$ of $P$ and every three distinct elements of $P$ are contained in at most one element of $B$;
 \item[{\rm (iii)}]
 Q is a subset of size $\delta$ of $P$ such that any triple of elements in $P$ with at most one element in $Q$, is contained in exactly one element of $B$;
 \end{itemize}
 
 {\bf Assumption} : Any such $\mathcal{D}$ is a $3 - (q^n + 1, q + 1, 1)$ design whenever $\delta \leq \delta_0$ with $\delta_0 \leq q - 2$.

 \begin{theorem}
 Consider a pseudo-hyperoval $\mathcal{O}$ in $\PG(3n - 1, q)$, $q = 2^h, h > 1$ and $n$ prime. Then $\mathcal{O}$ is regular if and only if at least $q^n + 1 - \delta_0$ elements of $\lbrace \Delta_0, \Delta_1, \cdots, \Delta_{q^n + 1} \rbrace$ are regular.
 \end{theorem}
{\em Proof}. \quad
Similar to the proof of Theorem 6.3. \eop

7.3. {\bf Acknowledgement} \\
We thank S. Rottey and G. Van de Voorde for several helpful discussions.

\end{document}